\documentclass[11pt,a4paper]{article}
\usepackage{amsmath}
\usepackage{amssymb}

\usepackage{graphicx}
\usepackage{amsfonts,amsmath, amssymb}

\numberwithin{equation}{section}

\usepackage{graphics}
\usepackage{epsfig}
\newtheorem{claim}{\bf \t}[part]

\newtheorem{theorem}{Theorem}[section]

\newtheorem{remark}[theorem]{Remark}

\def\t{\theta}

\begin{document}

\title{On Multi-Dimensional Sonic-Subsonic Flow }

\author{ {Feimin Huang,
Tianyi Wang, Yong Wang   }
\\
\ \\
   {\small \it Institute of Applied Mathematics,  Academy of Mathematics and Systems
   Science,}\\
{\small \it Academia Sinica, Beijing, 100190, P.R.China} \\
}
\date{ }
\maketitle

%\vskip 0.1cm \arraycolsep1.5pt
%\newtheorem{Lemma}{Lemma}[section]
%\newtheorem{Theorem}{Theorem}[section]
%\newtheorem{Definition}{Definition}[section]
%\newtheorem{Proposition}{Proposition}[section]
%\newtheorem{Remark}{Remark}[section]
%\newtheorem{Corollary}{Corollary}[section]
\begin{abstract}
In this paper, a compensated compactness framework is established
for sonic-subsonic approximate solutions to the
$n$-dimensional$(n\geq 2)$ Euler equations for steady irrotational
flow that may contain stagnation points. This compactness framework
holds provided that the approximate solutions are uniformly bounded
and satisfy $H^{-1}_{loc}(\Omega)$ compactness conditions. As
illustration, we show the existence of sonic-subsonic weak solution
to n-dimensional$(n\geq 2)$ Euler equations for steady irrotational
flow past obstacles or through an infinitely long nozzle. This is
the first result concerning the sonic-subsonic limit for
$n$-dimension$(n\geq 3)$.

\end{abstract}


\noindent\bf 2000 MR Subject Classification\quad \rm  35L65

%%%%%%%%%%%%%%%%%%%%%%%%%%%%%%%%%%%%%%%%%%%%%%%%%%%%%%%%%%%%%%%%%%%%%%%%%%%%%%%%%%%%%%%%%%%%%%%%%%
\section{Introduction}

The $n$-dimensional$(n\geq2)$ Euler equations for the steady
irrotational flow reads
\begin{eqnarray}\label{1.1}
\begin{cases}
\mbox{curl}\  u=0,\\
 \mbox{div}(\rho u)=0,\\
\mbox{div}(\rho u\otimes u+pI)=0,
\end{cases}
\end{eqnarray}
where $u=(u_1,\cdots,u_n)$  are the flow velocities,  and $\rho$ and
$p$ represent the density and pressure-density function,
respectively, and $(\mbox{curl}\
u)_{ij}=\partial_{x_j}u_i-\partial_{x_i}u_j,\ i,j=1,\cdots,n$, is a
$n\times n$ matrix and $I$  is a $n\times n$ unit matrix. Usually,
we require $p'(\rho)>0$ for $\rho>0$.

The famous Bernoulli's law can be easily derived from
$(\ref{1.1})_1$ and $(\ref{1.1})_3$:
\begin{equation}\label{1.5}
h(\rho)+\frac12q^2=const,
\end{equation}
where $q^2=|u|^2=\sum_{i=1}^{n} u_i^2$ is the flow speed and
$h(\rho)=\int_{1}^{\rho}\frac{p'(s)}{s}ds$ is the enthalpy.

In this paper, we are interested in the polytropic gas, that is
$p=p(\rho)=\frac{\rho^{\gamma}}{\gamma}$, for $\gamma
> 1$. Then $(\ref{1.5})$ is converted to the normalized formula
\begin{eqnarray}\label{1.2}
\rho=\rho(q)=(1-\frac{\gamma-1}{2}q^2)^{\frac{1}{\gamma-1}}.
\end{eqnarray}
 The local
sound speed is defined by
\begin{eqnarray}\label{1.3}
c^2=p'(\rho)=1-\frac{\gamma-1}{2}q^2.
\end{eqnarray}
At the sonic point $q=c$, \eqref{1.3} implies
$q^2=\frac{2}{\gamma+1}$. The critical speed $q_{cr}$ is defined as
$$q_{cr}=\sqrt{\frac{2}{\gamma+1}},$$
and the Bernoulli's law is rewritten as
$$q^2-q_{cr}^2=\frac{2}{\gamma+1}(q^2-c^2).$$
Thus the flow is subsonic when $q<q_{cr}$, sonic when $q=q_{cr}$,
and supersonic when $q>q_{cr}$.

For the isothermal flow, $p=\bar{c}^2\rho$, where $\bar{c}>0$ is the
constant sound speed, Bernoulli's law takes:
\begin{eqnarray}\label{1.4}
\rho=\rho(q)=\rho_0\exp{\bigg(-\frac{q^2}{2\bar{c}^2}\bigg)}
\end{eqnarray}
for some constant $\rho_0>0$.

It is well known that the steady irrotational Euler equations
\eqref{1.1} is of mixed type of partial differential equations,
which is elliptic if $q<q_{cr}$, parabolic if $q=q_{cr}$, hyperbolic
if $q>q_{cr}$.  Since the equations of uniform subsonic flow possess
ellipticity, its solutions have extra-smoothness to those related to
transonic flow or supersonic flow. There are a large of literatures
on the smooth uniform subsonic solutions, for instance, see
\cite{Bers2},\cite{Finn1},\cite{Frankl},\cite{Shiffman1},\cite{Shiffman2}
for two dimensional flow and
\cite{Dong1},\cite{Dong2},\cite{Xin3},\cite{Gilbarg1},\cite{Payne},\cite{Xin1},\cite{Xin2},\cite{Xin4}
for three dimensional flow. Among them, Frankl and Keldysh
\cite{Frankl} obtained the first result about the subsonic flow past
a two dimensional finite body (or airfoil). Shiffman in
\cite{Shiffman1}, \cite{Shiffman2} proved there exists a unique
subsonic potential flow around a given profile with finite energy
provided that the infinite free stream flow speed $q_\infty$ is less
than some critical speed, which was improved by Bers \cite{Bers2}.
Finn and Gilbarg \cite{Finn1} proved the uniqueness of the two
dimensional potential subsonic flow past a bounded obstacle with
given circulation and velocity at infinity. The first result for
three dimensional subsonic flow past an obstacle was given by Finn
and Gilbarg \cite{Gilbarg1} in which they studied the existence,
uniqueness and the asymptotic behavior with some restriction on Mach
number. Dong \cite{Dong1} extended the results of Finn and Gilbarg
\cite{Gilbarg1} to maximum Mach number $M < 1$ and to arbitrary
dimensions. Recently, Du, Xin and Yan \cite{Xin3} obtained the
smooth uniform subsonic for $n$-dimensional$(n\geq2)$ flow in an
infinitely long nozzle. For other related results, we refer to
\cite{Bers4,Bers5,Bers3,Bers1,Finn,Gilbarg2,Gilbarg3,Gilbarg4,Gilbarg5}
and references therein.

However, few result is known until now for the sonic-subsonic flow
and transonic flow, because the uniform ellipticity is lost and
shocks may present. That is, smooth solutions may not exist, and
weak solution is necessarily considered. Morawetz \cite{C.
Morawetz2,C. Morawetz3} firstly introduced the compensated
compactness method to study steady flow of irrotational Euler
equations. Indeed, Morawetz established a compactness framework
under assumption that the stagnation points and cavitation points
are excluded. Morawetz's result was improved by Chen, Slemrod and
Wang \cite{Chen8} in which the approximate solutions are constructed
by a viscous perturbation.

On the other hand, the first compactness framework on sonic-subsonic
irrotational flow in two dimension was recently due to Chen,
Dafermos, Slemrod and Wang \cite{Chen6} and Xie-Xin \cite{Xin1}
independently by combining the mass conservation, momentum, and
irrotational equations. The key point of \cite{Chen6} and
\cite{Xin1} is based on the fact that the two dimensional steady
flow can be regarded as the one dimensional system of conservation
laws, that is, $x$ is regarded as time $t$, so that the div-curl
lemma can be applied to the two momentum equations. In fact, the
authors \cite{Chen6} first applied the momentum equations to reduce
the support of the corresponding Young measure to two points, then
the irrotational equation and the mass equation are used to deduce
the Young measure to a Dirac measure.  Xie and Xin \cite{Xin2}
investigated the sonic-subsonic limit for the three-dimensional
axis-symmetric flow(it is similar to the two dimensional case)
through an axis-symmetric nozzle.

However, the compactness  framework established in \cite{Chen6} is
no longer effective for $n$-d $(n\geq 3)$ steady irrotational Euler
equations, which can not be reduced to one dimensional system of
conservation laws, and the famous div-curl lemma is no longer valid
for the momentum equations. In this paper, we find that it is
enough, by only using the mass conservation equation and
irrotational equations, to reduce the Young measure to a Dirac
measure for arbitrary dimension. Thus we establish a compactness
framework of approximate solutions for steady irrotational flow in
$n$-dimension$(n\geq2)$. It is worthy to point out that the famous
Bernoulli's law plays a key role in our proof. As application, we
show the sonic-subsonic limit
 for steady irrotational flows past
obstacles or through an infinitely long nozzle in
$n$-dimension$(n\geq2)$.

The rest of this paper is organized as follows. In section 2, we
establish the compactness framework of sonic-subsonic approximate
solutions for the system of  steady irrotational equations in
$n$-dimension$(n\geq2)$.  In section 3, we give two applications of
the compactness framework to show the existence of sonic-subsonic
flow over obstacles or through an infinitely long nozzle. Finally,
in section 4, we extend the compensated compactness framework to the
realm of self-similar solutions to the Euler equations for unsteady
irrotational flow.

\section{Compensated Compactness Framework for  Steady Irrotaitional  Flow in $n$-Dimension}

Let a sequence of function $u^\varepsilon(x)=(u^\varepsilon_1,
\cdots, u^\varepsilon_n)(x)$, defined on open subset $\Omega\subset
\mathbb{R}^n$, satisfy the following Conditions:

(A.1) $q^\varepsilon(x)=|u^\varepsilon(x)|\leq q_{cr}$ a.e. in
$\Omega$;

(A.2) $\mbox{curl}\  u^\varepsilon$, and
$\mbox{div}(\rho(q^\varepsilon) u^\varepsilon)$ are confined in a
compact set in $H_{loc}^{-1}(\Omega)$.

\

Based on the above conditions,  the famous div-curl lemma and the
Young measure representation theorem for a uniformly bounded
sequence of functions imply:
\begin{eqnarray}\label{2.1}
<\rho(q)q^2,
\nu(u)>=\sum\limits_{i=1}^n<\rho(q)u_i,\nu(u)><u_i,\nu(u)>,
\end{eqnarray}
where $\nu=\nu_{x}(u)$ is the associated Young measure (a
probability  measure) for the sequence
$u^\varepsilon(x)=(u^\varepsilon_1, \cdots, u^\varepsilon_n)(x)$.
Now, the main effort is to establish a compensated compactness
framework, namely, to prove that $\nu$ is a Diac measure by using
the identity $(\ref{2.1})$. This in turn implies the compactness of
the sequence $u^\varepsilon(x)=(u^\varepsilon_1, \cdots,
u^\varepsilon_n)(x)$ in $L^1_{loc}(\Omega)$.

\begin{theorem}(Compensated  compactness framework)\label{thm2.1}
Let a sequence of functions $u^\varepsilon(x)=(u^\varepsilon_1,
\cdots, u^\varepsilon_n)(x)$ satisfy conditions $(A.1)$ and $(A.2)$.
Then the associated Young measure $\nu$ is a Dirac mass and the
sequence $u^\varepsilon(x)$ is compact in $L^1_{loc}(\Omega)$; that
is, there is a subsequence (still labeled)
$u^\varepsilon(x)\rightarrow u(x)=(u_1, \cdots, u_n)(x)$ a.e. as
$\varepsilon\rightarrow 0$ and satisfying $q(x)=|u(x)|\leq q_{cr}$,
a.e. $x\in \Omega$.
\end{theorem}

\noindent\textbf{Proof}. Let
\begin{eqnarray}\label{2.3}
I(u^{(1)},u^{(2)})=\sum_{i=1}^n(u_i^{(1)}-u_i^{(2)})(\rho(q^{(1)})u_i^{(1)}-\rho(q^{(2)})u_i^{(2)}),
\end{eqnarray}
where $u^{(i)}=(u_1^{(i)}, \cdots, u_n^{(i)})$ and
$q^{(i)}=|u^{(i)}|$ for $i=1,2$ be two independent vector variables.

After some basic calculations on $I(u^{(1)},u^{(2)})$, we have
\begin{eqnarray}\label{2.4}
I(u^{(1)},u^{(2)})=\rho(q^{(1)})[(q^{(1)})^2-\sum\limits_{i=1}^n
u_i^{(1)}u_i^{(2)}]+\rho(q^{(2)})[(q^{(2)})^2-\sum\limits_{i=1}^n
u_i^{(1)}u_i^{(2)}].
\end{eqnarray}

Then the Cauchy inequality implies
\begin{eqnarray}\label{2.5}
I(u^{(1)},u^{(2)})&\geq
&\rho(q^{(1)})[(q^{(1)})^2-q^{(1)}q^{(2)}]+\rho(q^{(2)})[(q^{(2)})^2-q^{(1)}q^{(2)}]\nonumber\\
&=&(q^{(1)}-q^{(2)})(\rho(q^{(1)})q^{(1)}-\rho(q^{(2)})q^{(2)})\nonumber\\
&=&(q^{(1)}-q^{(2)})^2 \frac{d(\rho q)}{d q}(\tilde{q}).
\end{eqnarray}
where $\tilde{q}$ lies between $q^{(1)}$ and $q^{(2)}$ due to the
mean-value theorem. The famous Bernoulli's law \eqref{1.2} gives
that for $\gamma>1$,
$$
\rho(q)=(1-\frac{\gamma-1}{2}q^2)^{\frac{1}{\gamma-1}},
$$
which immediately implies
\begin{eqnarray}\label{2.6}
\frac{d(\rho q)}{d
q}=(1-\frac{\gamma-1}{2}q^2)^{\frac{1}{\gamma-1}-1}(1-\frac{\gamma+1}{2}q^2)=(1-\frac{\gamma-1}{2}q^2)^{\frac{1}{\gamma-1}-1}(1-\frac{q^2}{q_{cr}^2}).
\end{eqnarray}

For $\gamma = 1$, the Bernoulli's law is
$$
\rho(q)=\rho_0\exp{\bigg(-\frac{q^2}{2q_{cr}^2}\bigg)},
$$
which gives
\begin{eqnarray}\label{2.7}
\frac{d(\rho q)}{d
q}=\rho_0(1-\frac{q^2}{q_{cr}^2})\exp{\bigg(-\frac{q^2}{2q_{cr}^2}\bigg)}.
\end{eqnarray}

Thus , for sonic-subsonic flows, namely, $q^{(1)}, q^{(2)}\leq
q_{cr}$, $\eqref{2.5}-\eqref{2.7}$ imply
\begin{equation}\label{2.8}
I(u^{(1)},u^{(2)})=(q^{(1)}-q^{(2)})^2 \frac{d(\rho q)}{d
q}(\tilde{q})\geq 0
\end{equation}
and
\begin{equation}\label{2.9}
(q^{(1)}-q^{(2)})^2 \frac{d(\rho q)}{d q}(\tilde{q})=0,\ \mbox{if
and only if}\ q^{(1)}=q^{(2)}.
\end{equation}

From the identity $(\ref{2.1})$, noticing that the Young measure
$\nu$ is a probability measure, we have
\begin{eqnarray}\label{2.2}
<I(u^{(1)},u^{(2)}), \nu(u^{(1)})\otimes\nu(u^{(2)})>=0,
\end{eqnarray}
which together with \eqref{2.8} and \eqref{2.9} implies
$q^{(1)}=q^{(2)}$, where $\nu(u^{(1)})\otimes\nu(u^{(2)})$ is
understood as a product measure of $\nu(u^{(1)})$ and
$\nu(u^{(2)})$. With the property $q^{(1)}=q^{(2)}$ at hand, we have
from (\ref{2.3})
\begin{eqnarray}\label{2.10}
0&=&<I(u^{(1)},u^{(2)}),
\nu(u^{(1)})\otimes\nu(u^{(2)})>\nonumber\\
&=&<\rho(q^{(1)})\sum\limits_{i=1}^n
(u_i^{(1)}-u_i^{(2)})^2,\nu(u^{(1)})\otimes\nu(u^{(2)})>,
\end{eqnarray}
which immediately implies $u^{(1)}=u^{(2)}$, i.e, the Young measure
is a Dirac measure. This completes Theorem $\ref{thm2.1}$.

\

\begin{remark}\label{rem1} Theorem \ref{thm2.1} is valid
for any $n\geq2$. Namely, a compactness framework in Theorem
\ref{thm2.1} is established for sonic-subsonic limit for steady
irrotational flow in arbitrary dimension. From the Bernoulli's law
\eqref{1.5}, it is straightforward to extend Theorem \ref{thm2.1} to
the general pressure-density function $p$ satisfying $p'(\rho)>0$
for $\rho>0$.

Theorem \ref{thm2.1} is valid for any $n\geq2$. Namely, a
compactness framework in Theorem \ref{thm2.1} is established for
sonic-subsonic limit for steady irrotational flow in arbitrary
dimension. From the Bernoulli's law \eqref{1.5}, it is
straightforward to extend Theorem \ref{thm2.1} to the general
pressure-density function $p$ satisfying $p'(\rho)>0$ for $\rho>0$.
Futhermore, The above idea can be applied for the sonic-subsonic
limit of the self-similarity solutions, that is,
$(\rho,u)(x,t)=(\rho,u)(\frac{x}{t}), x\in R^n$, for unsteady
compressible Euler system. The idea is also available for
Euler-Poisson systems.

\end{remark}

\

We now consider a sequence of approximate solutions $u^\varepsilon$
to the Euler equations $(\ref{1.1})_1$, $(\ref{1.1})_2$ and the
Bernoulli's law $(\ref{1.2})$ or $(\ref{1.4})$. That is, besides
Conditions $(A.1)$ and $(A.2)$, the approximate solutions
$u^\varepsilon$ further satisfy
\begin{eqnarray}\label{2.11}
\begin{cases}
\mbox{curl}\  u^\varepsilon=o_1(\varepsilon),\\
 \mbox{div}(\rho(q^\varepsilon) u^\varepsilon)=o_2(\varepsilon),\\
\end{cases}
\end{eqnarray}
where $o_1(\varepsilon),o_2(\varepsilon)\rightarrow 0$ in the sense
of distributions as $\varepsilon\rightarrow 0$. Then, as a corollary
of Theorem \ref{thm2.1}, we have

\begin{theorem}(Convergence of approximate solutions)\label{thm2.2}
Let $u^\varepsilon(x)=(u^\varepsilon_1, \cdots, u^\varepsilon_n)(x)$
be a sequence of approximate solutions satisfying $(\ref{2.11})$ and
the Bernoulli's law $(\ref{1.2})$ or $(\ref{1.4})$. Then, there is a
subsequence (still labeled) $u^\varepsilon(x)$ that converges a.e.
as $\varepsilon\rightarrow 0$ to a weak solution $u(x)=(u_1, \cdots,
u_n)(x)$ to the Euler equations of $(\ref{1.1})_1$, $(\ref{1.1})_2$
and the Bernoulli's law $(\ref{1.2})$ or $(\ref{1.4})$ satisfying
$q(x)=|u(x)|\leq q_{cr}$, a.e. $x\in \Omega$.
\end{theorem}
\begin{remark}\label{rem2} For any functions
$Q(u)=(Q_1(u),\cdots,Q_n(u))$ satisfying
\begin{equation}\label{2.12}
\mbox{div}(Q(u^\varepsilon))=o(\varepsilon),
\end{equation}
where $o(\varepsilon)\rightarrow 0$ in the sense of distributions as
$\varepsilon\rightarrow 0$, from the strong convergence of
$u^\varepsilon$, $Q(u)=0$ holds in the sense of distributions. So if
we have
\begin{equation}
\mbox{div}(\rho(q^\varepsilon) u^\varepsilon\otimes
u^\varepsilon+p^\varepsilon I)=o(\varepsilon)\rightarrow 0,\
\mbox{in the sense of distributions},
\end{equation}
the weak solution in Theorem $\ref{thm2.2}$ also satisfies the
momentum equations $(\ref{1.1})_3$ in the sense of distributions.
\end{remark}

\

There are various ways to construct approximate solutions by either
numerical methods or analytical methods such as vanishing viscosity
methods. In the next section, we will show two examples to apply the
compactness framework built in Theorem \ref{thm2.1}.

\section{Sonic Limit of Irrotational Subsonic Flows in n-Dimension}

In this section, we wish to apply the compactness framework
established in Theorem \ref{thm2.1} to obtain the sonic limit of
$n$-dimensional$(n\geq2)$ steady irrotational subsonic flows.

\begin{figure}[htbp]
\small \centering
\includegraphics[width=10cm]{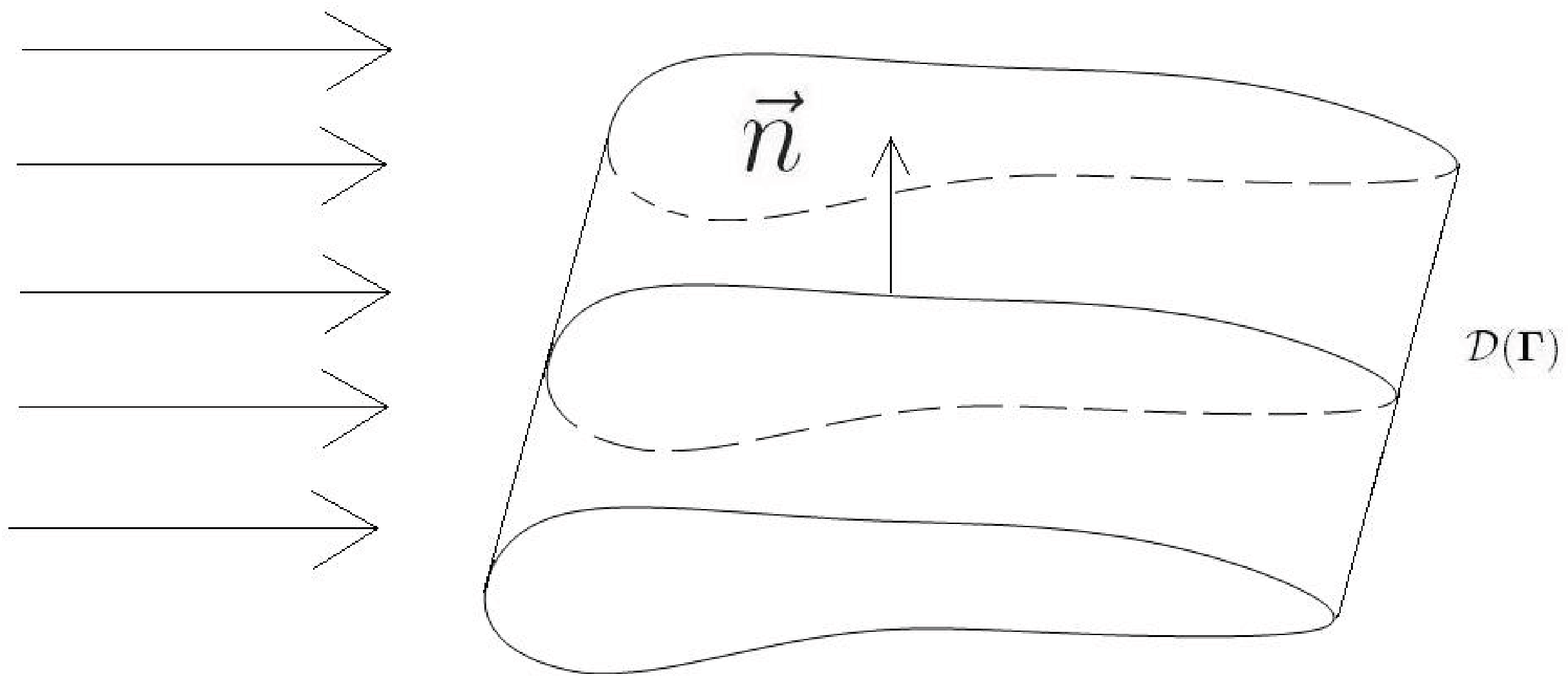}
\caption{General high dimensional case.}
 \label{Fig2}
\end{figure}

Firstly, we give an example of subsonic-sonic limit past obstacles.
Let the obstacle $\Gamma$ be one or several closed $n-1(n\geq2)$
dimensional hypersurfaces. We shall always assume $\Gamma\in
C^{2,\tau_0}$. Denote by $\mathcal{D}(\Gamma)$ the domain exterior
to $\Gamma$, see  Fig \ref{Fig2},

\

 $\mbox{Problem}\ \mathbf{\mathbf{P}}(u_\infty):$  Let $n\geq2$. Find functions
$u=(u_1,\cdots,u_n)$ satisfy $(\ref{1.1})_1,(\ref{1.1})_2$ with the
Bernoulli's law \eqref{1.2} or \eqref{1.4}, and the slip boundary
condition
\begin{equation}\label{3.1}
(\rho u)\cdot \vec{n}=0\ \ \mbox{on}\ \ \Gamma,
\end{equation}
where $\vec{n}$ denotes the unit inward normal of domain
$\mathcal{D}(\Gamma)$, and the limit
\begin{equation}\label{3.2}
u_\infty=\lim_{|x|\rightarrow \infty}u(x),
\end{equation}
exists and is finite.

\

The main result of \cite{Dong1} is described as follows:

\begin{theorem}[Uniform Subsonic Flows Past An Obstacle for  n-D Case
\cite{Dong1}]\label{thm3.2}
 Let \\ $q_{\infty}:=|u_\infty|$. There
exists a positive number $\hat{q}<q_{cr}$, so that
$\mathbf{\mathbf{P}}(u_\infty)$ has a uniform subsonic solutions if
$0\leq q_\infty<\hat{q}$. Furthermore, let
$q_m(q_\infty)=\sup\limits_{x\in\mathcal{D}(\Gamma)}|u(x)|$, then
the function $q_m(q_\infty)\in C[0,\hat{q})$ and
$q_m(q_\infty)\rightarrow q_{cr}$ as $q_\infty\rightarrow \hat{q}$.
\end{theorem}

We note that the above Theorem does not apply to the critical flows,
that is, those flows for which $q_m(q_\infty)= q_{cr}$ and which
hence must be sonic at some point. Here, by Theorem $\ref{thm2.1}$,
we establish a more general result.

\begin{theorem}[Sonic Limit Past An Obstacle]\label{thm3.4}
Let $u_{\infty}^\varepsilon\rightarrow \hat{u}_\infty\ \mbox{as}\
\varepsilon\rightarrow0$ be a sequence of speeds at $\infty$ with
$q_\infty^{\varepsilon}<\hat{q}=|\hat{u}_\infty|$, and
$u^{\varepsilon}=(u_1^{\varepsilon},\cdots,u_n^{\varepsilon})$ be
the corresponding solutions to  Problem
$\mathbf{P}(u_\infty^\varepsilon)$. Then, as
$u_{\infty}^\varepsilon\rightarrow \hat{u}_\infty$, the solution
sequence $u^{\varepsilon}(x)$ possess a subsequence (still denoted
by) $u^{\varepsilon}(x)$ that converge a.e. in $\mathcal{D}(\Gamma)$
to a vector function
 $u(x)=(u_1,\cdots,u_n)(x)$ which is a weak solution of Problem $\mathbf{P}(u_\infty)$ with $q_\infty=\hat{q}$.
Furthermore the limit velocity $u=(u_1,\cdots,u_n)$ also satisfies
$(\ref{1.1})_3$ in the sense of
 distributions and the
boundary conditions \eqref{3.1} as the normal trace of the
divergence-measure field $(u_1,\cdots,u_n)$ on the boundary (see
\cite{Chen7}).
\end{theorem}
\noindent{\bf Proof}. The strong solutions $u^\varepsilon$ satisfy
\eqref{1.1}, and the Bernoulli's law and are uniform subsonic
solutions of Problem $\mathbf{P}(u_\infty^\varepsilon)$. Hence
Theorem $\ref{thm2.1}$ immediately implies that the Young measure is
a Dirac mass and the convergence is strong a.e. in
$\mathcal{D}(\Gamma)$. The boundary conditions are satisfied for $u$
in the sense of Chen-Frid \cite{Chen7}. On the other hand, Since
$(\ref{1.1})_3$ holds for the sequence of subsonic solutions
$u^{\varepsilon}(x)$, it is straightforward to see that $u$ also
satisfies $(\ref{1.1})_3$ in the sense of
 distributions. This completes the proof of Theorem \ref{thm3.4}.

\

\begin{figure}[htbp]
\small \centering
\includegraphics[width=10cm]{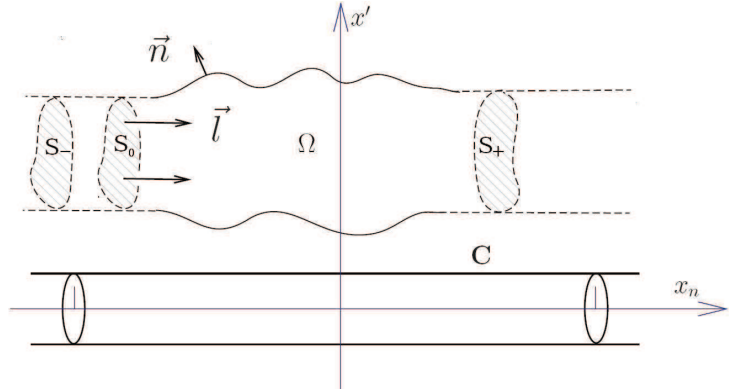}
\caption{n-dimensional nozzle} \label{Fig3}
\end{figure}

Now we give another example of subsonic-sonic limit through an
infinite long nozzle. As in \cite{Xin3}, denote the
multi-dimensional nozzle domain by $\Omega$ which satisfies the
following regularity assumption: there exists an invertible
$C^{2,\alpha}$ map $T: \bar{\Omega}\rightarrow\bar{C}:x\rightarrow
y$ satisfying
\begin{eqnarray}\label{3.3}
\begin{cases}
T(\partial\Omega)=\partial C,\\
\mbox{For any}\ k\in\mathbb{R},\
T(\Omega\cap\{x_n=k\})=B(0,1)\times\{y_n=k\},\\
\|T\|_{C^{2,\alpha}}, \|T^{-1}\|_{C^{2,\alpha}}\leq K<\infty,\\
\mbox{The nozzle approaches to a cylinder in the far fields},i.e,\\
\Omega\cap\{x_n=k\}\rightarrow S_{\pm}\ \mbox{as}\ k\rightarrow
\pm\infty,\mbox{respectively},
\end{cases}
\end{eqnarray}
where $K$ is a uniform constant, $C=B(0,1)\times(-\infty,+\infty)$
is a unit cylinder in $\mathbb{R}^n$, $B(0,1)$ is unit ball in
$\mathbb{R}^{n-1}$ centered at the origin, $S_{\pm}$ are $n-1$
dimensional simply connected $C^{2,\alpha}$, $x_n$ is the
longitudinal coordinate and $x'=(x_1,\cdots,x_{n-1})\in
\mathbb{R}^{n-1}$, see Fig \ref{Fig3},

\

$\mbox{Problem}\ \mathbf{\mathbf{\tilde{P}}}(m_0):$  Let $n\geq2$.
Find functions $u=(u_1,\cdots,u_n)$ satisfy
$(\ref{1.1})_1,(\ref{1.1})_2$ with the Bernoulli's law \eqref{1.2}
or \eqref{1.4}, and the slip boundary condition
\begin{equation}\label{3.4}
(\rho u)\cdot \vec{n}=0\ \ \mbox{on}\ \ \partial\Omega,
\end{equation}
where $\vec{n}$ denotes the unit outward normal of domain $\Omega$;
and the mass flux condition
\begin{equation}\label{3.5}
\int_{S_0}\rho(|u|^2)u\cdot \vec{l}dS=m_0>0
\end{equation}

\

The main result of \cite{Xin3} is stated as follows:
\begin{theorem}[Uniform Subsonic Flows in n-D Nozzle\cite{Xin3}]\label{thm3.3}
There is a critical mass flux\\  $M_c>0$, which depends only on
$\Omega$, such that if $0<m_0<M_c$, then
$\mathbf{\mathbf{\tilde{P}}}(m_0)$ has a unique uniformly subsonic
flow through the nozzle, i.e,
$q_m(m_0):=\sup\limits_{x\in\Omega}|u(x)|<q_{cr}$. The velocity
$(u_1,\cdots,u_n)$ is Holder continuous. Moreover,
$q_m(m_0)\rightarrow q_{cr}$ as $m_0\rightarrow M_c$.
\end{theorem}

\

Similar to Theorem \ref{thm3.4}, we have
\begin{theorem}[Sonic Limit Through A Nozzle]\label{thm3.5}
Let $0<m_0^{\varepsilon}<M_c$ be a sequence of mass flux, and let
 $u^\varepsilon(x)$ be the corresponding solution to $\mathbf{\mathbf{\tilde{P}}}(m_0^\varepsilon)$. Then as
 $m_0^{\varepsilon}\rightarrow M_c$, the solution
sequence $u^{\varepsilon}(x)$ possess a subsequence (still denoted
by) $u^{\varepsilon}(x)$ that converge strongly a.e. in $\Omega$ to
a vector function
 $u(x)=(u_1,\cdots,u_n)(x)$ which is a weak solution of   $\mathbf{\mathbf{\tilde{P}}}(M_c)$  with
Bernoulli's law. Furthermore the limit velocity $u=(u_1,\cdots,u_n)$
also satisfies $(\ref{1.1})_3$ in the sense of
 distributions and the boundary conditions \eqref{3.4}
as the normal trace of the divergence-measure field
$(u_1,\cdots,u_n)$ on the boundary (see \cite{Chen7} ).

\end{theorem}

\begin{remark}  In this section we only give two examples as
applications of compactness framework established in Theorem
\ref{thm2.1}. Certainly it can be used to other cases as long as
conditions $(A.1)$ and $(A.2)$ are satisfied.
\end{remark}

\

\noindent {\bf Acknowledgments:} The first author would like to
thank Prof. Cathleen S.Morawetz for her kind hospitality and helpful
discussion when he visited Courant Institute in 2006. The research
of FMH was supported in part by NSFC Grant No. 10825102 for
distinguished youth scholar, and National Basic Research Program of
China (973 Program) under Grant No. 2011CB808002.

\

\end{document}